\documentclass[10pt]{amsart}

\usepackage{amsmath,epsfig}
\usepackage{amssymb}
\usepackage{amsthm}
\usepackage{amsfonts}
\usepackage[all]{xy}

\newcommand{\mgn}{\overline{M}_{g,P}}

\newcommand{\mgbar}{\overline{\mathcal{M}}_{g,P}}

\newcommand{\mgpbar}{\overline{\mathcal{M}}_{g, P\cup\{q\}}}

\newtheorem{theorem}{Theorem}
\newtheorem{corollary}[theorem]{Corollary}

\newtheorem{remark}[theorem]{Remark}

\title{Chern Classes of the Moduli Stack of Curves}\author{Gilberto Bini}
\email{gilberto.bini@mat.unimi.it} \curraddr{Dipartimento di
Matematica, Universit\`a degli Studi di Milano \\ Via C. Saldini 50 \\
20133 Milano \\ Italy.}

\begin{document}
\begin{abstract}
Here we calculate the Chern classes of ${\overline{\mathcal
M}}_{g,n}$, the moduli stack of stable $n$-pointed genus $g$ curves. In
particular, we prove that such classes lie in the tautological ring.
\end{abstract}
\maketitle

\section{Introduction}\label{sec:0}

Let $g$ and $n$ be non-negative integers such that $n > 2-2g$. We
denote by ${\overline{\mathcal M}}_{g,n}$ the Deligne-Mumford stack of
stable genus $g$ curves with $n$ marked points. More generally, if $P$
is a set with $n$ elements, it will be technically convenient to work
with $\mgbar$, i.e., the stack of genus $g$ stable curves whose
marked points are labelled by $P$. The natural projection $\varphi$
from the stack $\mgbar$ to the coarse moduli space $\mgn$ induces an
isomorphism $\varphi_*$ at the level of Chow rings - hereafter we
shall only deal with rational coefficients. Following \cite{FP}, we
denote by $R^*(\mgn)$ the tautological ring of $\mgn$. By abuse of
notation, we shall denote the image of $R^*(\mgn)$ under
$(\varphi_*)^{-1}$ by the same symbol.

Tautological classes have been intensely studied in the last few
years: see \cite{rav} for a synoptic survey on the most recent
developments in this area. In particular, it is not at all clear which
classes lie in the tautological ring. In fact, constructions of
tautological or non-tautological classes can be very diverse in
nature: combinatorial - as conjectured by Kontsevich and proved in,
e.g., \cite{mon} - or purely algebro-geometric (e.g.,
\cite{GraPan:01}).  As shown by Mumford in \cite{Mum}, the
Grothendieck-Riemman-Roch Theorem allows one to express Chow classes
in terms of tautological ones. Incidentally, this is done for the
canonical class of $\mgbar$ in \cite{HarM}. That calculation can be
rephrased in terms of stacks. For foundational material on stacks we
refer the reader to \cite{vis} and, especially for Chern classes, to
\cite{man}.

In the present paper, we extend Mumford's work and calculate all Chern
classes of the (smooth) moduli stack of curves, i.e., of the tangent
bundle ${\mathcal T}_{\mgbar}$ to $\mgbar$. In spite of their
geometric significance, these classes have not been hitherto
computed. As above, we apply the Grothendieck-Riemann-Roch Theorem. In
fact, we refine combinatorial arguments, and we manage to get explicit
formulas. This shows in particular that such classes are tautological,
thus yielding new elements in $R^*(\mgn)$.

Our presentation is rather concise since most of the theoretic
material has been explored by several authors. We briefly recall the
needed preliminaries in Section \ref{sec:1} and prove the main result
in Section 3. Finally, we give some examples and compare our formulas
with previous results.

Throughout, we shall work over the field of complex numbers.

\section{Preliminaries}\label{sec:1}

Let $\mgbar$ be the moduli stack of $P$-pointed stable curves of genus
$g$. As usual, $\pi: {\mathcal C} \rightarrow \mgbar$ will denote the
universal curve with sections $\sigma_p$ for $p \in P$. Set, further,
$D:=\sum_{p \in P} \sigma_{p_*}(1)$. The universal cotangent classes
on $\mgbar$ are defined as $\psi _{p}=c_{1}(\sigma
_{p}^{*}(\omega_{\pi }))$ for $p\in P$, where $\omega_{\pi}$ is the
relative dualizing sheaf of $\pi$. The collection of all moduli stacks
$\mgbar$ is equipped with some natural morphisms, namely:
\begin{equation}
\label{csi}
\xi _{G }:\prod _{v\in V}\overline{\mathcal {M}}_{g(v),l(v)}\rightarrow \overline{\mathcal {M}}_{g,P},
\end{equation}
where $G$ is a stable graph - see, e.g., \cite{AC2}. For example, it
is well known that the boundary $\partial {\overline{\mathcal
M}}_{g,P}$ can be described in terms of the following morphisms:
\begin{equation}
\label{bounda2}
\xi_{irr}: {\overline{\mathcal M}}_{g-1,P\cup\{q_1,q_2\}} \rightarrow {\overline{\mathcal M}}_{g,P},
\end{equation}
\begin{equation} \xi_{h,A}: {\overline{\mathcal M}}_{h,A \cup \{r_1\}} 
\times {\overline{\mathcal M}}_{g-h,A^c \cup \{r_2\}} \rightarrow {\overline{\mathcal M}}_{g,P},
\label{bounda1}
\end{equation}
where $0\leq h \leq g$, $A \subseteq P$, and both $2h-1+|A|$ and
$2(g-h)-1+|A^{c}|$ are positive. Finally, let $\delta$ be the boundary
class defined as
\begin{equation}
\label{totbou}
\delta=\frac{1}{2}\xi_{irr \, \, *}(1)+\frac{1}{2}\sum _{h}\sum
_{A\subseteq P}\xi _{h,A \, \, *}(1).
\end{equation}

Let $K=c_{1}\left(\omega_{\pi }(D)\right)$. Following \cite{AC1}, the
Mumford classes on $\mgbar$ are defined as $ \kappa _{m}=\pi
_{*}(K^{m+1})$. For $P =\emptyset$ their analogue was first introduced
by Mumford in \cite{Mum}. Another generalization of Mumford's
$\kappa_m$'s to the case of $P$-pointed curves is given by the classes
$ \widetilde{\kappa }_{m}=\pi _{*}(c_{1}(\omega _{\pi })^{m+1})$. As
proved in \cite{AC1}, the following relationship holds:
\begin{equation}
\label{mor1}
\kappa_{m}=\widetilde{\kappa }_{m}+\sum _{p\in P}\psi _{p}^{m}.
\end{equation}

The Hodge bundle ${\mathbb E}$ on $\mgbar$ is defined as
$\pi_*\omega_{\pi}$. From \cite{Bin} and \cite{Mum} we have
\begin{equation*}
ch({\mathbb E})= g + \frac12 \sum_{m \geq 1}\frac{B_{2m}}{(2m)!}
                   \Bigl\{\widetilde{\kappa}_{2m-1} +
\end{equation*}
\begin{equation*}
\xi_{irr \, \, *}(\psi_{q_1}^{2m-2}-\psi_{q_1}^{2m-3}\psi_{q_2} +
\ldots +\psi_{q_2}^{2m-2})+
\end{equation*} 
\begin{equation*}
\sum_{h=0}^{g}\sum_{A\subseteq P}\xi_{G_{h,A}\, \, *} (\psi_{r_1}^{
2m-2}\otimes 1 -\psi_{r_1}^{2m-3} \otimes \psi_{r_2} + \ldots
+1\otimes \psi_{r_2}^{2m-2})\Bigr\},
\end{equation*}
where $B_{2m}$ are the Bernoulli numbers\footnote{We recall that Bernoulli numbers are defined as follows:
$$
\frac{x}{e^x-1}=1-\frac{1}{2}x + \sum_{m \geq 1} \frac{B_{2m}}{(2m)!}x^{2m}.
$$}.
Note in particular that
\begin{equation}
\label{lambda}
ch_1({\mathbb E})=c_1({\mathbb E}):=\lambda= \frac{1}{12}
\left(\kappa_1 - \sum_{p \in P}\psi_p\right) +\delta.
\end{equation} 

As first shown by Mumford in \cite{Mum}, the Grothendieck-Riemann-Roch
Theorem (the G-R-R Theorem for short) can be applied to the universal
curve $\pi: {\mathcal C} \rightarrow \mgbar$. Alternatively, one can
use the G-R-R Theorem stated in \cite{Toe}. For the sake of
completeness, we report the statement below.

\begin{theorem}(G-R-R Theorem)
\label{crucial} Let ${\mathcal E}$ be a locally free sheaf on ${\mathcal C}$. 
Then 
$$ ch(\pi_!{\mathcal E})=\pi_*\left(ch({\mathcal
E})Td^{\vee}(\Omega_{\pi})\right),
$$ where $\Omega_{\pi}$ is the sheaf of relative 
K\"{a}hler differentials.
\end{theorem}

In Proposition \ref{crucial}, $ch({\mathcal E})$ and $Td({\mathcal
E})$ denote the Chern character and the dual Todd class of ${\mathcal E}$,
respectively. Some formulas for these classes can be found, for
instance, in \cite{ACGH}. Here, we just remark two basic facts. First,
notice that
\begin{equation}
\label{link} 
ch({\mathcal E}^*)= rk({\mathcal E})+\sum_{j \geq 1}
(-1)^jch_j({\mathcal E}).
\end{equation} 

Second, let $\mu=(1^{m_1} 2^{m_2} \ldots
i^{m_i}\ldots)$ be a partition of weight $j$, where $j$ is a positive
integer. Define $ch_{\mu}({\mathcal E})$ to be the product $
ch_1^{m_1}({\mathcal E})ch_2^{m_2}({\mathcal E}) \ldots
ch_i^{m_i}({\mathcal E})\ldots $ As proved in \cite{mac}, (2.14'), the
following holds:
\begin{equation}
\label{eccoti}
c_j({\mathcal E})= \sum_{\mu \vdash j}(-1)^{j-l(\mu)} \prod_{r \geq 1} \frac{((r-1)!)^{m_r}}{m_r!}ch_{\mu}({\mathcal E}), \quad j \geq 1,
\end{equation}
where the sum ranges over all partitions $\mu$ of $j$, and $l(\mu)$ is
the length of $\mu$.

Finally, we recall some properties of $Z$, the singular locus of
$\pi$. For more details the reader is referred to \cite{FP2}. $Z$ is a
closed substack of codimension $2$ in ${\mathcal C}$. Moreover, there
exists a double \`{e}tale covering $\varepsilon: \widetilde{Z}
\rightarrow Z$ obtained from the choice of the branches incident at
the nodes corresponding to points in $Z$. Let $\iota: \widetilde{Z}
\rightarrow {\mathcal C}$ be the natural composition. Denote by
${\mathcal L}$ and ${\mathcal L}'$ the line bundles corresponding to
the cotangent directions along the branches. Then $\varepsilon^*({\mathcal
N}_{Z})= {\mathcal L}\oplus {\mathcal L}'$, where ${\mathcal N}_Z$ is
the normal bundle of $Z$ in ${\mathcal C}$. Moreover, note that $\pi
\circ \iota$ maps $\widetilde{Z}$ onto $\partial \mgbar$. In other
words, we have
\begin{equation}
\label{whoknows}
\pi \circ \iota = \xi_{irr} + \sum_h 
\sum_{\substack{A \subseteq P}} \xi_{h,A},
\end{equation}
where both $2h-1+|A|$ and $2(g-h)-1+|A^{c}|$ are positive.

\section{The Chern Character of $\mgbar$}\label{sec2}

In this section we apply the G-R-R Theorem to the sheaf
$\Omega_{\pi}(D)\otimes \omega_{\pi}$. By standard duality theorems
and deformation theory \cite{HarM},
${\pi}_*\left(\Omega_{\pi}(D)\otimes \omega_{\pi}\right)$ is the
cotangent bundle on $\mgbar$. In the sequel, we closely follow
\cite{Mum}, p. 302 ff.

We recall that
\begin{equation}
\label{tweede}
ch\left(\Omega_{\pi}(D)\right)= ch\left(\omega_{\pi}(D)\right) - 
ch\left({\mathcal O}_Z\right),
\end{equation}
and
\begin{equation}
\label{second}
Td^{\vee}(\Omega_{\pi})=
Td^{\vee}(\omega_{\pi})\left(Td^{\vee}({\mathcal O}_Z)\right)^{-1},
\end{equation}
where ${\mathcal O}_Z$ is viewed as a sheaf on ${\mathcal C}$. For the
purpose of what follows, we need to determine the power series
$ch({\mathcal O}_Z)Td^{\vee}({\mathcal O}_Z)^{-1}$. Since
$Td^{\vee}({\mathcal O}_Z)^{-1}$ is a polynomial in the Chern characters 
$ch_k({\mathcal O}_Z)$, the G-R-R Theorem applied to $\iota$ yields a
universal power series $\Theta$. For all $\nu: Y \rightarrow X$, an
inclusion of a smooth codimension two subvariety in a smooth variety,
the series $\Theta$ satisfies the following identity:
\begin{equation}
\label{codim} ch({\mathcal O}_Y)Td^{\vee}({\mathcal O}_Y)^{-1}=
\nu_*\left(\Theta\left(c_1\left({\mathcal N}_Y\right), c_2\left({\mathcal N}_Y\right)\right)\right).
\end{equation}

To compute $\Theta$, take $Y=D_1D_2$. Then the exact sequence
$$ 
0 \rightarrow {\mathcal O}_X(-D_1 -D_2) \rightarrow {\mathcal
O}_X(-D_1) \oplus {\mathcal O}_X(-D_2) \rightarrow {\mathcal O}_X
\rightarrow {\mathcal O}_Y \rightarrow 0
$$
yields
\begin{equation}
\label{thisisnew}
ch({\mathcal O}_Y)= (1-e^{-D_1})(1-e^{-D_2}).
\end{equation} 

By \cite{Mum}, p. 303, we have
\begin{equation}
\label{feree}
Td^{\vee}({\mathcal O}_Y)^{-1}= \frac{D_1 D_2}{D_1+D_2} \frac{1 -
e^{-D_1-D_2}}{(1-e^{-D_1})(1-e^{-D_2})}.
\end{equation}

Therefore, we get
\begin{equation}
ch({\mathcal O}_Y)Td^{\vee}({\mathcal O}_Y)^{-1}=D_1D_2 \sum_{j \geq
1} (-1)^{j-1}\frac{(D_1+D_2)^{j-1}}{j!}.
\end{equation}

Thus, the following holds:
\begin{equation}
\label{sithistime}
\Theta(D_1+D_2, D_1D_2)= \sum_{j \geq 1} (-1)^{j-1}\frac{(D_1+D_2)^{j-1}}{j!}.
\end{equation}

As noted in \cite{FP2}, \eqref{thisisnew} can be applied to the
morphism $\iota: \widetilde{Z} \rightarrow {\mathcal C}$ as well. In
this case, we need to replace the Chern classes of ${\mathcal N}_Z$
with those of $\varepsilon^*{\mathcal N}_Z$.

We now state the main result of this section.

\begin{theorem}
\label{main}
The Chern character of the cotangent bundle on $\mgbar$ is given by
\begin{eqnarray}
ch\left({\mathcal T}_{\mgbar}^*\right) &=& \displaystyle{\sum_{j \geq
1} \frac{\kappa_{j-1}}{j!}} + \frac12\displaystyle{\sum_{t \geq
1}\frac{\kappa_t}{t!}} - \displaystyle{\sum_{m \geq 3}a_m
\kappa_{m-1}} + \\ \nonumber \\ \nonumber & & ch({\mathbb E}) -1 -\frac12
\xi_{irr \, \, *}\left(\Xi^{(1)}\right) - \frac12\sum_{h=0}^{g}\sum_{A\subseteq
P}\xi_{G_{h,A}\, \, *}\left(\Xi^{(2)}\right),\\ \nonumber
\end{eqnarray}
where
\begin{equation}
\label{aemme}
a_m=\sum_{h=1}^{\lfloor \frac{m-1}{2} \rfloor}
\frac{B_{2h}}{(2h)!(m-2h)!}\,,
\end{equation}

\begin{equation}
\label{prima}
\Xi^{(1)} =\sum_{k \geq 1}(-1)^{k-1} \frac{\left(\psi_{q_1}+
\psi_{q_2}\right)^{k-1}}{k!}
\end{equation}
and
\begin{equation}
\label{seconda}
\Xi^{(2)} = \sum_{k \geq 1}(-1)^{k-1} \frac{\left(\psi_{r_1} \otimes 1 +
1 \otimes \psi_{r_2}\right)^{k-1}}{k!}
\end{equation}  
\end{theorem}

{\bf Proof}. By \eqref{tweede} and \eqref{second}, we have
\begin{eqnarray}
\label{first}
ch\left(\bigl(\Omega_{\pi}(D)\otimes
\omega_{\pi}\bigr)Td^{\vee}(\Omega_{\pi})\right) &=&
ch(\omega_{\pi})Td^{\vee}(\omega_{\pi})Td^{\vee}({\mathcal
O}_Z)^{-1}ch(\omega_{\pi}(D)) \\ &-& 
ch(\omega_{\pi})Td^{\vee}(\omega_{\pi})Td^{\vee}({\mathcal
O}_Z)^{-1}ch({\mathcal O}_Z).\nonumber
\end{eqnarray}

As proved in \cite{AC1}, the first Chern class of $\omega_{\pi}(D)$ is
equal to $\psi_q$. Here $\psi_q$ denotes the universal cotangent class
on $\mgpbar$ corresponding to the marked point $q$.  \newpage

Since ${\mathcal
O}_Z$ is supported on $Z$ and the marked points are non-singular, we
get

\begin{eqnarray}
\label{fourteen}
ch(\omega_{\pi})Td^{\vee}(\omega_{\pi})Td^{\vee}({\mathcal
O}_Z)^{-1}ch(\omega_{\pi}(D)) &=&
ch(\omega_{\pi})Td^{\vee}(\omega_{\pi})ch(\omega_{\pi}(D)) \\ \nonumber &+&
Td^{\vee}({\mathcal O}_Z)^{-1}-1,
\end{eqnarray}
and
\begin{equation}
\label{fijvtien}
ch(\omega_{\pi})Td^{\vee}(\omega_{\pi})Td^{\vee}({\mathcal
O}_Z)^{-1}ch({\mathcal O}_Z)= Td^{\vee}({\mathcal
O}_Z)^{-1}ch({\mathcal O}_Z).
\end{equation}

In particular, it is easy to check that \eqref{fourteen} is equal to
\begin{equation}
ch(\omega_{\pi})Td^{\vee}(\omega_{\pi})\left[ch(\omega_{\pi}(D))-1\right] +ch(\omega_{\pi})Td^{\vee}(\Omega_{\pi}).
\label{comoantes}
\end{equation}

By definition, we have

\begin{equation}
\label{fifteen}
\begin{array}{c}
ch(\omega_{\pi}(D)) = e^{\psi_q}, \\ \\
ch(\omega_{\pi})Td^{\vee}(\omega_{\pi}) = e^{(\psi_q-D)}
\displaystyle{\frac{\psi_q-D}{e^{(\psi_q-D)}-1}}=
\displaystyle{\frac{D-\psi_q}{e^{D-\psi_q}-1}}.
\end{array}
\end{equation}

We recall that the formal expansion of the second power series in
\eqref{fifteen} is given by

\begin{equation}
\label{sixteen}
1 + \frac12(\psi_q-D)+\displaystyle{\sum_{j \geq
1}\frac{B_{2j}}{(2j)!}(-\psi_q+D)^{2j}},
\end{equation}
where $B_{2j}$ is the $(2j)$-th Bernoulli number.  Since $\psi_q\cdot D=0$,
the first term of \eqref{comoantes} is given by

\begin{equation}
\label{eighteen}
\sum_{j \geq 1} \frac{\psi_q^j}{j!} + \frac12\displaystyle{\sum_{t \geq 1}\frac{\psi_q^{t+1}}{t!}} - \displaystyle{\sum_{m \geq 3}a_m \psi_q^m}\,,
\end{equation}
where
$a_m$ is defined in \eqref{aemme}.

Let us apply ${\pi}_*$ to \eqref{first}. The contribution in
\eqref{comoantes} yields

\begin{equation}
\label{twenty}
 \displaystyle{\sum_{j \geq 1} \frac{\kappa_{j-1}}{j!}} +
\frac12\displaystyle{\sum_{t \geq 1}\frac{\kappa_t}{t!}} -
\displaystyle{\sum_{m \geq 3}a_m \kappa_{m-1}} + ch({\mathbb E}) -1.
\end{equation}

On the other hand, the contribution in \eqref{fijvtien} is given by
\begin{equation}
\label{twentythree}
\frac12(\pi \circ \iota)_* \left(\Theta\left(c_1\left(\varepsilon^*{\mathcal
N}_Z\right), c_2\left(\varepsilon^*{\mathcal N}_Z\right)\right)\right).
\end{equation}

By \eqref{whoknows}, this is equivalent to
$$ \frac12 \xi_{irr\, \,
*}\left(\Xi^{(1)}\right)+ \frac12\sum_{h=0}^{g}\sum_{A\subseteq
P}\xi_{G_{h,A}\, \, *}\left(\Xi^{(2)}\right),
$$ where $\Xi^{(1)}$ and $\Xi^{(2)}$ are defined in \eqref{prima} and
\eqref{seconda}, respectively.

\begin{flushright}
$\Box$
\end{flushright}

\begin{remark}
\label{nota}
By Theorem \ref{main} and \eqref{lambda}, we get
$$ 
ch_0\left({\mathcal T}^*_{\mgbar}\right)= rk({\mathcal T}^*_{\mgbar})=3g-3+n,
$$
$$ ch_1\left({\mathcal T}^*_{\mgbar}\right)= K_{\mgbar}=13\lambda +
\sum_{p\in P} \psi_p -2\delta,
$$ where $K_{\mgbar}$ is the canonical class of the stack $\mgbar$ and
$\delta$ is defined in \eqref{totbou}.  These formulas agree with
previous known results: see, e.g., \cite{Bin2}, \cite{HarM}.
\end{remark}

\begin{remark}
Note that Theorem \ref{main} and Formula \eqref{link} give the Chern
character of ${\mathcal T}_{\mgbar}$.
\end{remark}

By Formula \eqref{eccoti}, we obtain an expression for the Chern classes
of $\mgbar$.
\begin{corollary}
\label{lomando}
For $j\geq 1$, we have 
\begin{equation}
c_j(\mgbar)= \sum_{\mu \vdash j}(-1)^{j-l(\mu)}
\prod_{r \geq 1} \frac{((r-1)!)^{m_r}}{m_r!}ch_{\mu}(\mgbar).
\end{equation}
\end{corollary}

\begin{flushright}
$\Box$
\end{flushright}

{\bf Example.} We give some examples for low $j$'s. In higher degrees,
one can use John Stembridge's symmetric function package SF for {\tt
maple} \cite{ste}. For the sake of simplicity, we denote $\sum_p
\psi_p$ by $\psi$. As noted in Remark \ref{nota}, we have
$$
c_1(\mgbar)= -13\lambda - \psi +2\delta.
$$
From Corollary \ref{lomando}, we get
$$
ch_2(\mgbar)= \frac{\kappa_2}{3} + \frac14\xi_{irr \, \,
*}\left(\psi_{q_1} + \psi_{q_2}\right) + \frac14\sum_{h,A}\xi_{h,A \, \,
*}\left(\psi_{r_1}\otimes 1 + 1 \otimes \psi_{r_2}\right).
$$
hence
\begin{eqnarray*}
c_2(\mgbar) &=& \frac12 (-13\lambda - \psi +2\delta)^2- \frac13
\kappa_2 - \frac 14 \xi_{irr \, \, *}\left(\psi_{q_1} +
\psi_{q_2}\right) - \\ & & \frac14\sum_{h,A}\xi_{h,A \, \,
*}\left(\psi_{r_1}\otimes 1 + 1 \otimes \psi_{r_2}\right).
\end{eqnarray*}

Finally, the degree $3$ Chern character is equal to
\begin{eqnarray*}
ch_3(\mgbar)&=& - \frac{\kappa_3}{12} - ch_3({\mathbb E}) +
\frac{1}{12}\xi_{irr \, \, *}\left(\psi_{q_1}^2 + \psi_{q_1}\psi_{q_2}
+ \psi_{q_2}^2\right) + \\ \\ && \frac{1}{12}\sum_{h,A}\xi_{h,A \, \,
*}\left(\psi_{r_1}^2\otimes 1 + \left(\psi_{r_1}\otimes
\psi_{r_2}\right) + 1 \otimes \psi_{r_2}^2\right).
\end{eqnarray*}

Thus, we get
\begin{eqnarray*}
c_3(\mgbar)&=& \left[+13\lambda + \psi
-2\delta\right]\left[\frac13
\kappa_2 + \frac14\xi_{irr \, \,
*}\left(\psi_{q_1} + \psi_{q_2}\right)\right. \\ & & \left.  + \frac14\sum_{h,A}\xi_{h,A \, \,
*}\left(\psi_{r_1}\otimes 1 + 1 \otimes \psi_{r_2}\right)\right] - \\ \\ & & 
\frac16\left[-13\lambda - \psi +2\delta\right]^3 + 2ch_3(\mgbar). 
\end{eqnarray*}

As a result of Theorem \ref{main} and the
definition of tautological classes, we get new elements in the
tautological ring of $\mgbar$. Precisely, the following holds.

\begin{corollary}
The Chern classes of $\mgbar$ are tautological.
\end{corollary}

\begin{flushright}
$\Box$
\end{flushright}


\end{document}